\documentclass{article}
\usepackage{amssymb,latexsym,amsmath,amsthm,mathrsfs}
\usepackage{graphicx}
\newcommand{\comment}[1]{}
\begin{document}
\title{On the sums of series of reciprocals\footnote{Presented to the
St. Petersburg Academy on December 5, 1735.
Originally published as
{\em De summis serierum reciprocarum}, Commentarii academiae scientiarum Petropolitanae \textbf{7} (1740), 123--134. E41 in the Enestr{\"o}m index. Translated from the Latin by Jordan Bell, School of Mathematics and Statistics, Carleton University, Ottawa, Canada. Email: jbell3@connect.carleton.ca}}
\author{Leonhard Euler}
\date{}
\maketitle

\S 1. So much work has been done on the series of the reciprocals of powers of the natural numbers, that it seems hardly likely to be able to discover anything new about them. For nearly all
who have thought about the sums of series have also inquired into the sums of this kind of series, and yet have not been able by any means to express them in
a 
convenient form.
I too, in spite of repeated efforts, in which I attempted various methods for summing these series, could achieve nothing more than approximate values for their sums or reduce them to the quadrature of highly transcendental curves; 
the former of these is described in the next article, and the latter fact I have set out in preceding ones. 
I speak here about the series of fractions whose numerators are 1, and indeed whose denominators are the squares, or the cubes, or other ranks, of the natural numbers; of this type are $1+\frac{1}{4}+\frac{1}{9}+\frac{1}{16}+\frac{1}{25}+$ etc., likewise $1+\frac{1}{8}+\frac{1}{27}+\frac{1}{64}+$ etc. and similarly for higher powers, whose general terms are contained in the form $\frac{1}{x^n}$.

\S 2. I have recently found, quite unexpectedly, an elegant expression for the sum of this series $1+\frac{1}{4}+\frac{1}{9}+\frac{1}{16}+$ etc., which depends on the quadrature of the circle, so that if the true sum of this series is obtained, from it at once the quadrature of the circle follows. Namely, I have found for six times the sum of this series to be equal to the square of the perimeter of a circle whose diameter is 1; or by putting the sum of this series $=s$, then $\surd 6s$ will hold to 1 the ratio of the perimeter to the diameter.
Indeed I recently showed for the sum of this series to be approximately $1,644934066842264364$; multiplying this number by six, and then taking the square root, the very same number $3,141592653589793238$ is found which expresses the perimeter of a circle whose diameter is 1.
Following again the same steps by which I had arrived at this sum, I have discovered that the sum of the series $1+\frac{1}{16}+\frac{1}{81}+\frac{1}{256}+\frac{1}{625}+$ etc. also depends on the quadrature of the circle. 
Namely, the sum of this multiplied by 90 gives the biquadrate of the perimeter of a circle whose diameter is 1.
And by similar reasoning I have likewise been able to determine the sums of the subsequent series in which the exponents are even numbers.

\begin{figure}
\begin{center}
\includegraphics[scale=0.75]{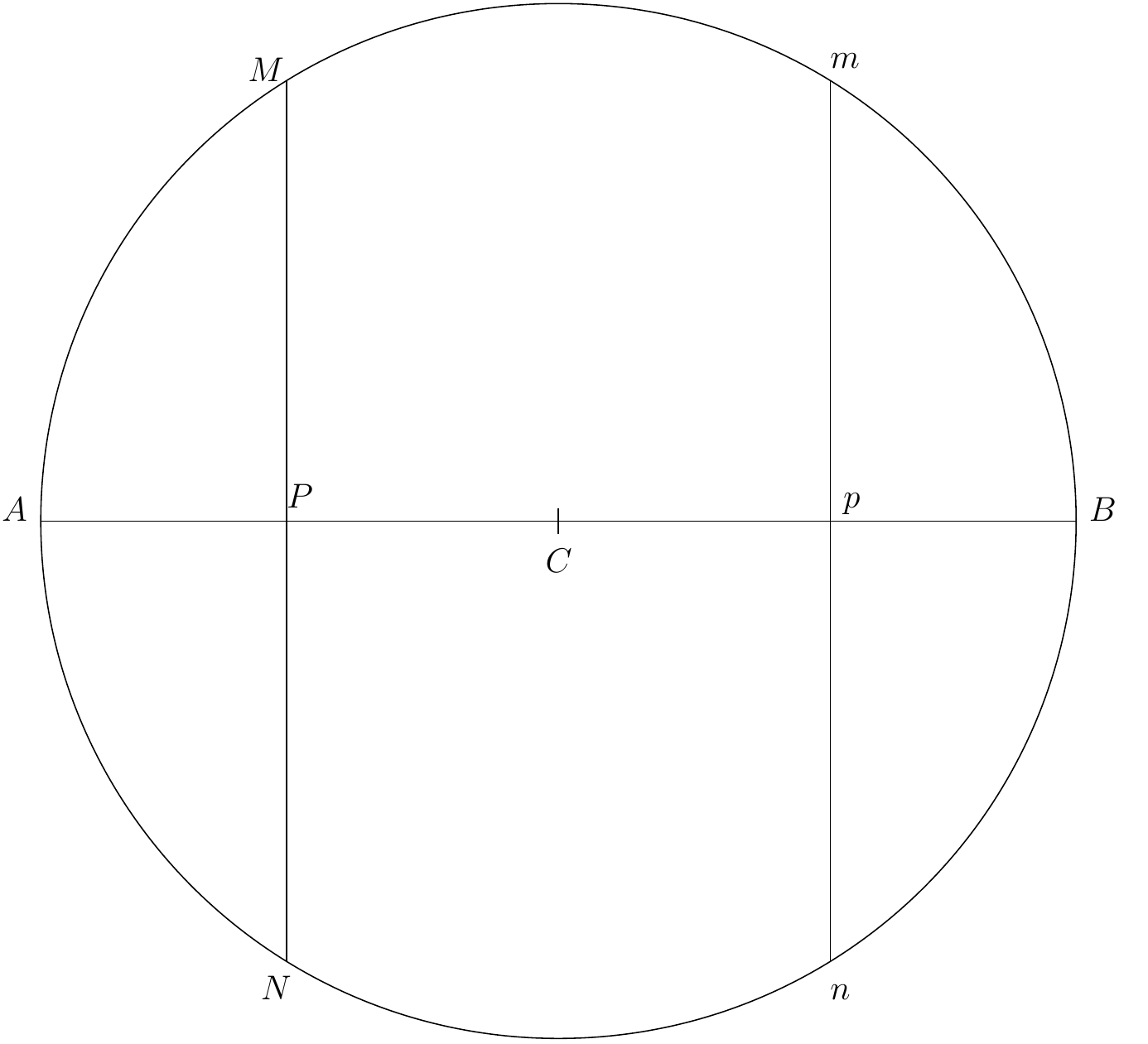}
\end{center}
\end{figure}

\S 3. I will therefore explain most thoroughly how I accomplished this;  
I will set forth everything in the order in which I have made use of it.
In the circle $AMBNA$ described with center $C$, radius $AC$ or $BC=1$, I have considered an arc $AM$, whose sine is $MP$ and whose cosine indeed is $CP$. Now, by putting the arc $AM=s$, sine $PM=y$, and cosine $CP=x$, then by a well known method the sine $y$, and likewise the cosine $x$, can be defined by a series from the given arc $s$, as it can be seen that in general $y=s-\frac{s^3}{1\cdot 2\cdot 3}+\frac{s^5}{1\cdot 2\cdot 3\cdot 4\cdot 5}-\frac{s^7}{1\cdot 2\cdot 3\cdot 4\cdot 5\cdot 6\cdot 7}+$ etc. and also $x=1-\frac{s^2}{1\cdot 2}+\frac{s^4}{1\cdot 2\cdot 3\cdot 4}-\frac{s^6}{1\cdot 2\cdot 3\cdot 4\cdot 5\cdot 6}+$ etc. 
It was by considering these very equations that I was led to the above
sums of the series of reciprocals in question;
and in fact, either of these equations can equally well be used to get to the same end, so it therefore suffices to handle just one of them, which I shall set forth.

\S 4. The first equation $y=s-\frac{s^3}{1\cdot 2\cdot 3}+\frac{s^5}{1\cdot 2\cdot 3\cdot 4\cdot 5}-\frac{s^7}{1\cdot 2\cdot 3\cdot 4\cdot 5\cdot 6\cdot 7}+$ etc. thus expresses a relation between the arc and sine. For from it,
from a given arc its sine can be determined, and likewise for a given sine the arc. But now I will consider the sine $y$ as given, and I will investigate how the arc must arise. Indeed, it should first be noted that since
innumerable arcs correspond to the same sine $y$, the given equation will satisfy innumerable arcs. Certainly, if in this equation $s$ is seen as an unknown, it will have infinitely many degrees, and it will not be surprising that this equation should contain countless simple factors which, when put equal to zero, then yield a suitable value for $s$.

\S 5. Furthermore, if all the factors of this equation were known, then all the roots of it, or values of $s$, would be known, and on the other hand if all the values of $s$ are assigned then too all the factors of this equation will be possessed. Thus, since it is simpler for me to determine the factors than the roots, I transform the given equation into this form:
$0=1-\frac{s}{y}+\frac{s^3}{1\cdot 2\cdot 3\cdot y}-\frac{s^5}{1\cdot 2\cdot 3\cdot 4\cdot 5\cdot y}+$ etc. If now all the roots of this equation, or all the arcs with this same sine $y$, were $A,B,C,D,E$, etc., then likewise the factors will be all these quantities, $1-\frac{s}{A},1-\frac{s}{B},1-\frac{s}{C},1-\frac{s}{D}$, etc. Wherefore it will be $1-\frac{s}{y}+\frac{s^3}{1\cdot 2\cdot 3\cdot y}-\frac{s^5}{1\cdot 2\cdot 3\cdot 4\cdot 5\cdot y}+\textrm{etc.}=(1-\frac{s}{A})(1-\frac{s}{B})(1-\frac{s}{C})(1-\frac{s}{D})$ etc.

\S 6. But from the very nature of solving equations, it is apparent that the coefficient of the term in which $s$ appears, or $\frac{1}{y}$, is equal to the sum of all the coefficients of $s$ in the factors, that is, $\frac{1}{y}=\frac{1}{A}+\frac{1}{B}+\frac{1}{C}+\frac{1}{D}+$ etc.
Next,
the coefficient of $s^2$, which is $=0$ because this term does not appear in the equation, is equal to the sum of the factors from two terms of the series, $\frac{1}{A},\frac{1}{B},\frac{1}{C},\frac{1}{D}$, etc. In turn,
$-\frac{1}{1\cdot 2\cdot 3\cdot y}$ will be equal to the sum of the factors from three terms of the series $\frac{1}{A},\frac{1}{B},\frac{1}{C},\frac{1}{D}$, etc. And in the same manner, it will be that $0=$ the sum of the factors from four terms of the series, and $+\frac{1}{1\cdot 2\cdot 3\cdot 4\cdot 5\cdot y}=$ the sum of the factors from five terms of this series, and so on.

\S 7. But by letting $AM=A$ be the least arc whose sine is $PM=y$, and the semiperimeter of the circle $=p$, then $A,p-A,2p+A,3p-A,4p+A,5p-A,6p+A$, etc. and similarly $-p-A,-2p+A,-3p-A,-4p+A,-5p-A$, etc. will be all the arcs with the same sine $y$. Since we have supposed before that the series $\frac{1}{A},\frac{1}{B},\frac{1}{C},\frac{1}{D}$, etc. is transformed into this $\frac{1}{A}$, $\frac{1}{p-A}$, $\frac{1}{-p-A}$, $\frac{1}{2p+A}$, $\frac{1}{-2p+A}$, $\frac{1}{3p-A}$, $\frac{1}{-3p-A}$, $\frac{1}{4p+A}$, $\frac{1}{-4p+A}$, etc., the sum of all these terms is therefore $=\frac{1}{y}$; furthermore, the sum of the factors from two terms of this series is equal to $0$; the sum of the factors from three $=\frac{-1}{1\cdot 2\cdot 3\cdot y}$; the sum of the factors from four $=0$; the sum of the factors from five $=\frac{+1}{1\cdot 2\cdot 3\cdot 4\cdot 5\cdot y}$; the sum of the factors from six $=0$. And so on.

\S 8. On the other hand, if an arbitrary series $a+b+c+d+e+f+$ etc. is considered, the sum of which is $\alpha$, the sum of the factors from two terms is $=\beta$, the sum of the factors
from three terms $=\gamma$, the sum of the factors from four terms $=\delta$, etc., then the sum of the squares of all the terms will be $a^2+b^2+c^2+d^2+\textrm{etc.}=\alpha^2-2\beta$; the sum of the cubes $a^3+b^3+c^3+d^3+\textrm{etc.}=\alpha^3-3\alpha \beta+3\gamma$; the sum of the biquadrates $=\alpha^4-4\alpha^2 \beta+4\alpha \gamma+2\beta^2-4\delta$. To make it clearer how these formulae proceed, we put the sum of these terms $a,b,c,d$, etc. to be $=P$, the sum of the squares of the terms $=Q$, the sum of the cubes $=R$, the sum of the biquadrates $=S$, the sum of the fifth powers $=T$, the sum of the sixth powers $=V$, etc. With this done, it will be $P=\alpha$; $Q=P\alpha -2\beta$; $R=Q\alpha-P\beta +3\gamma$; $S=R\alpha-Q\beta+P\gamma-4\delta$; $T=S\alpha-R\beta+Q\gamma-P\delta+5\epsilon$; etc.

\S 9. As in our case of the series $\frac{1}{A},\frac{1}{p-A},\frac{1}{-p-A},\frac{1}{2p+A},\frac{1}{-2p+A},\frac{1}{3p-A},\frac{1}{-3p-A}$, etc., the sum of all the terms, or $\alpha$, would be $=\frac{1}{y}$; the sum of the factors from two or $\beta=0$, and further for higher numbers of terms, $\gamma=\frac{-1}{1\cdot 2\cdot 3\cdot y}$; $\delta=0$; $\epsilon=\frac{+1}{1\cdot 2\cdot 3\cdot 4\cdot 5\cdot y}$; $\zeta=0$; etc., therefore the sum of all the terms in this series will be $P=\frac{1}{y}$; the sum of all the squares of these terms $Q=\frac{P}{y}=\frac{1}{y^2}$; the sum of all the cubes of these terms $R=\frac{Q}{y}-\frac{1}{1\cdot 2\cdot y}$; the sum of all the biquadrates $S=\frac{R}{y}-\frac{P}{1\cdot 2\cdot 3\cdot y}$. And then in turn $T=\frac{S}{y}-\frac{Q}{1\cdot 2\cdot 3\cdot y}+\frac{1}{1\cdot 2\cdot 3\cdot 4\cdot y}$; $V=\frac{T}{y}-\frac{R}{1\cdot 2\cdot 3\cdot y}+\frac{P}{1\cdot 2\cdot 3\cdot 4\cdot 5\cdot y}$; $W=\frac{V}{y}-\frac{S}{1\cdot 2\cdot 3\cdot y}+\frac{Q}{1\cdot 2\cdot 3\cdot 4\cdot 5\cdot y}-\frac{1}{1\cdot 2\cdot 3\cdot 4\cdot 5\cdot 6\cdot y}$. From this rule, the sums of the remaining higher powers are easily determined.

\S 10. We shall now put the sine $PM=y$ equal to the radius, so that $y=1$. The least arc $A$ whose sine is $1$ will be a quarter part of the perimeter, $=\frac{1}{2}p$, or by denoting with $q$ a quarter part of the perimeter, it will be $A=q$ and $p=2q$. The above series therefore turns into this, $\frac{1}{q}$, $\frac{1}{q}$, $-\frac{1}{3q}$, $-\frac{1}{3q}$, $+\frac{1}{5q}$, $+\frac{1}{5q}$, $-\frac{1}{7q}$, $-\frac{1}{7q}$, $+\frac{1}{9q}$, $+\frac{1}{9q}$, etc., where the terms arise as equal pairs. Thus the sum of these terms, which is $\frac{2}{q}(1-\frac{1}{3}+\frac{1}{5}-\frac{1}{7}+\frac{1}{9}-\frac{1}{11}+\textrm{etc.})$, is equal to $P=1$ itself. It therefore follows from this that $1-\frac{1}{3}+\frac{1}{5}-\frac{1}{7}+\frac{1}{9}-\frac{1}{11}+\textrm{etc.}=\frac{q}{2}=\frac{p}{4}$. Thus four times this series is equal to the semiperimeter of a circle whose radius is $1$, or the total perimeter of a circle whose diameter is $1$. And indeed this is the very same series discovered some time ago by {\em Leibniz}, by which he defined the quadrature of the circle. From this,  
if our method should appear to some as not reliable enough,
a great confirmation comes to light here; thus there should not be any doubt about the rest that will be derived from this method.

\S 11. We shall now take up the squares of the terms that arise in the case when $y=1$, and the series $+\frac{1}{q^2}+\frac{1}{q^2}+\frac{1}{9q^2}+\frac{1}{9q^2}+\frac{1}{25q^2}+\frac{1}{25q^2}+$ etc. will appear, whose sum is $\frac{1}{q^2}(\frac{1}{1}+\frac{1}{9}+\frac{1}{25}+\frac{1}{49}+\textrm{etc.})$, which will then be equal to $Q=P=1$. From this it follows that the sum of the series $1+\frac{1}{9}+\frac{1}{25}+\frac{1}{49}+$ etc. is $=\frac{q^2}{2}=\frac{p^2}{8}$, where $p$ denotes the total perimeter of a circle whose diameter is $=1$. But the sum of this series $1+\frac{1}{9}+\frac{1}{25}+$ etc. determines the sum of the series
$1+\frac{1}{4}+\frac{1}{9}+\frac{1}{16}+\frac{1}{25}+$ etc., because
the latter series minus a quarter of itself gives the former, and
therefore, the sum of the former series plus a third of itself is equal to the sum of the latter series.
Wherefore, it will be $1+\frac{1}{4}+\frac{1}{9}+\frac{1}{16}+\frac{1}{25}+\frac{1}{49}+\textrm{etc.}=\frac{p^2}{6}$, and indeed, the sum of this series multiplied by $6$ is equal to the square of the perimeter of a circle whose diameter is $1$; this is the very proposition that I had mentioned at the start.

\S 12. Since therefore in the case when $y=1$ it holds that $P=1$ and $Q=1$, the values of the other letters will be as follows: $R=\frac{1}{2}$; $S=\frac{1}{3}$; $T=\frac{5}{24}$; $V=\frac{2}{15}$; $W=\frac{61}{720}$; $X=\frac{17}{335}$; etc. Then furthermore, since the sum of the cubes is equal to $R=\frac{1}{2}$, it will be $\frac{2}{q^3}(1-\frac{1}{3^3}+\frac{1}{5^3}-\frac{1}{7^3}+\frac{1}{9^3}-\textrm{etc.})=\frac{1}{2}$. Whence it will be $1-\frac{1}{3^3}+\frac{1}{5^3}-\frac{1}{7^3}+\frac{1}{9^3}-\textrm{etc.}=\frac{q^3}{4}=\frac{p^3}{32}$. Therefore, the sum of this series multiplied by $32$ gives the cube of the perimeter of a circle whose diameter is $1$. Similarly, the sum of the biquadrates, which is $\frac{2}{q^4}(1+\frac{1}{3^4}+\frac{1}{5^4}+\frac{1}{7^4}+\frac{1}{9^4}+\textrm{etc.})$ should be equal to $\frac{1}{3}$, and thus it will be $1+\frac{1}{3^4}+\frac{1}{5^4}+\frac{1}{7^4}+\frac{1}{9^4}+\textrm{etc.}=\frac{q^4}{6}=\frac{p^4}{96}$. But in fact, this series multiplied by $\frac{16}{15}$ is equal to the very series $1+\frac{1}{2^4}+\frac{1}{3^4}+\frac{1}{4^4}+\frac{1}{5^4}+\frac{1}{6^4}+$ etc., whence that series is equal to $\frac{p^4}{90}$; that is to say, the sum of the series $1+\frac{1}{2^4}+\frac{1}{3^4}+\frac{1}{4^4}+$ etc. multiplied by 90 gives the biquadrate of the perimeter of a circle whose diameter is 1.

\S 13. The sums of the above powers can be found in a similar way; it will in fact turn out to follow that $1-\frac{1}{3^5}+\frac{1}{5^5}-\frac{1}{7^5}+\frac{1}{9^5}-\textrm{etc.}=\frac{5q^5}{48}=\frac{5p^5}{1536}$; and also $1+\frac{1}{3^6}+\frac{1}{5^6}+\frac{1}{7^6}+\frac{1}{9^6}+\textrm{etc.}=\frac{q^6}{15}=\frac{p^6}{960}$. And indeed having found the sum of this series, we can immediately perceive the sum of the series $1+\frac{1}{2^6}+\frac{1}{3^6}+\frac{1}{4^6}+\frac{1}{5^6}+$ etc., which will be $=\frac{p^6}{945}$. Next for the seventh powers it will be $1-\frac{1}{3^7}+\frac{1}{5^7}-\frac{1}{7^7}+\frac{1}{9^7}-\textrm{etc.}=\frac{61q^7}{1440}=\frac{61p^7}{184320}$, and for the eighth, $1+\frac{1}{3^8}+\frac{1}{5^8}+\frac{1}{7^8}+\frac{1}{9^8}+\textrm{etc.}=\frac{17q^8}{630}=\frac{17p^8}{161280}$; whence it is deduced that $1+\frac{1}{2^8}+\frac{1}{3^8}+\frac{1}{4^8}+\frac{1}{5^8}+\frac{1}{6^8}+\textrm{etc.}=\frac{p^8}{9450}$. It should be observed, however, that in these series of powers, the signs of the terms of odd exponents alternate, while for the even powers they are in fact equal; because of this, the sum of the general series $1+\frac{1}{2^n}+\frac{1}{3^n}+\frac{1}{4^n}+$ etc. can be exhibited only for the cases when $n$ is an even number.
Besides this, it should be noted that when the general term of the series $1,1,\frac{1}{2},\frac{1}{3},\frac{5}{24},\frac{2}{15},\frac{61}{720},\frac{17}{315}$, etc. can be determined, whose values we have found as the letters $P,Q,R,S$, etc., then with it the quadrature of the circle can be given.

\S 14. In this we have so far put the sine $PM$ equal to the radius. We can also see what sort of series is produced if other values are taken for $y$. Thus, if $y=\frac{1}{\surd 2}$, the least arc corresponding to this sine is $\frac{1}{4}p$. Therefore by putting $A=\frac{1}{4}p$, the series of simple terms, or first powers, will be $\frac{4}{p}+\frac{4}{3p}-\frac{4}{5p}-\frac{4}{7p}+\frac{4}{9p}+\frac{4}{11p}-$ etc. The sum $P$ of this series is equal to $\frac{1}{y}=\surd 2$. It will therefore be obtained that $\frac{p}{2\surd 2}=1+\frac{1}{3}-\frac{1}{5}-\frac{1}{7}+\frac{1}{9}+\frac{1}{11}-\frac{1}{13}-\frac{1}{15}+$ etc.; the rule of the signs of this series differs from that of {\em Leibniz}, and in fact has been published before by {\em Newton}. Indeed the sum of the squares of these terms,  
namely $\frac{16}{p^2}(1+\frac{1}{9}+\frac{1}{25}+\frac{1}{49}+\textrm{etc.})$, is equal to $Q=2$. It will therefore be $1+\frac{1}{9}+\frac{1}{25}+\frac{1}{49}+\textrm{etc.}=\frac{p^2}{8}$, as was found earlier.

\S 15. If it were taken $y=\frac{\surd 3}{2}$, then the least arc corresponding to this sine will be $60^{\circ}$, and then $A=\frac{1}{3}p$. In this case the following series of terms will appear, $\frac{3}{p}+\frac{3}{2p}-\frac{3}{4p}-\frac{3}{5p}+\frac{3}{7p}+\frac{3}{8p}-$ etc., and the sum of these terms is equal to $\frac{1}{y}=\frac{2}{\surd 3}$ itself. Therefore we will have $\frac{2p}{3\surd 3}=1+\frac{1}{2}-\frac{1}{4}-\frac{1}{5}+\frac{1}{7}+\frac{1}{8}-\frac{1}{10}-\frac{1}{11}+$ etc. Indeed, the sum of the squares of these terms is $=\frac{1}{y^2}=\frac{4}{3}$; whence it follows that $\frac{4p^2}{27}=1+\frac{1}{4}+\frac{1}{16}+\frac{1}{25}+\frac{1}{49}+\frac{1}{64}+$ etc., where in this 
series the third terms are always missing.
And in fact, this series depends on the one $1+\frac{1}{4}+\frac{1}{9}+\frac{1}{16}+$ etc., whose sum was earlier found to be $=\frac{p^2}{16}$; for that series less a ninth part is equal to the above series, whose sum should therefore be $=\frac{p^2}{6}(1-\frac{1}{9})=\frac{4pp}{27}$. Similarly, if other sines are assumed, then other series arise, for the simple terms, the squares of the terms and too for higher powers, which likewise involve the quadrature of the circle.

\S 16. But if it were put $y=0$, then a series of this type could no longer be assigned because $y$ has been put in the denominator, that is, the original equation has been divided by $y$.
On the other hand though, other series can indeed be worked out in this situation, which however are series of the form $1+\frac{1}{2^n}+\frac{1}{3^n}+\frac{1}{4^n}+$ etc. if $n$ is an even number:
for the sums of these series which are to be found, I can already deduce them separately from this case where $y=0$. Now, if we let $y=0$, the fundamental equation becomes $0=s-\frac{s^3}{1\cdot 2\cdot 3\cdot y}+\frac{s^5}{1\cdot 2\cdot 3\cdot 4\cdot 5}-\frac{s^7}{1\cdot 2\cdot 3\cdot 4\cdot 5\cdot 6\cdot 7}+$ etc., and the roots of this equation give all the arcs whose sine is $=0$.
There is one minimal root $s=0$, so that when the equation is divided by $s$ it exhibits all the remaining arcs whose sine is $=0$, so that these arcs are all roots of the equation $0=1-\frac{s^2}{1\cdot 2\cdot 3}+\frac{s^4}{1\cdot 2\cdot 3\cdot 4\cdot 5}-\frac{s^6}{1\cdot 2\cdot 3\cdot 4\cdot 5\cdot 6\cdot 7}+$ etc. 
Of course these arcs whose sine is $=0$ are $p,-p,+2p,-2p,3p,-3p$, etc., of which the second of each pair is negative, which the equation itself also tells us, because the dimensions of $s$ are even.
Hence the divisors of this equation will be $1-\frac{s}{p},1+\frac{s}{p},1-\frac{s}{2p},1+\frac{s}{2p}$, etc., and by combining each of these pairs of divisors together it will be $1-\frac{s^2}{1\cdot 2\cdot 3}+\frac{s^4}{1\cdot 2\cdot 3\cdot 4\cdot 5}-\frac{s^6}{1\cdot 2\cdot 3\cdot 4\cdot 5\cdot 6\cdot 7}+\textrm{etc.}=(1-\frac{s^2}{p^2})(1-\frac{s^2}{4p^2})(1-\frac{s^2}{9p^2})(1-\frac{s^2}{16p^2})$ etc.

\S 17. Now it is evident from the nature of equations that the coefficient of $ss$, or $\frac{1}{1\cdot 2\cdot 3}$, is equal to $\frac{1}{p^2}+\frac{1}{4p^2}+\frac{1}{9p^2}+\frac{1}{16p^2}+$ etc. And indeed the sum of the factors from two terms of this series will be $=\frac{1}{1\cdot 2\cdot 3\cdot 4\cdot 5}$; and the sum of the factors from three $=\frac{1}{1\cdot 2\cdot 3\cdot 4\cdot 5\cdot 6\cdot 7}$, etc. Because of this, it will be, like in \S 8, $\alpha=\frac{1}{1\cdot 2\cdot 3}$; $\beta=\frac{1}{1\cdot 2\cdot 3\cdot 4\cdot 5}$; $\gamma=\frac{1}{1\cdot 2\cdot 3\cdot 4\cdot 5\cdot 6\cdot 7}$; etc. By also putting the sum of the terms $\frac{1}{p^2}+\frac{1}{4p^2}+\frac{1}{9p^2}+\frac{1}{16p^2}+\textrm{etc.}=P$, and the sum of the squares of these terms $=Q$; the sum of the cubes $=R$; the sum of the biquadrates $=S$; etc., it will be by \S 8, $P=\alpha=\frac{1}{1\cdot 2\cdot 3}=\frac{1}{6}$; $Q=P\alpha-2\beta=\frac{1}{90}$; $R=Q\alpha-P\beta+3\gamma=\frac{1}{945}$; $S=R\alpha-Q\beta+P\gamma-4\delta=\frac{1}{9450}$; $T=S\alpha-R\beta+Q\gamma-P\delta+5\epsilon=\frac{1}{93555}$; $V=T\alpha-S\beta+R\gamma-Q\delta+P\epsilon-6\zeta=\frac{691}{6825\cdot 93555}$, etc.

\S 18. From these observations the following sums are thus derived:
\[
\begin{split}
&1+\frac{1}{2^2}+\frac{1}{3^2}+\frac{1}{4^2}+\frac{1}{5^2}+\textrm{etc.}=\frac{p^2}{6}=P\\
&1+\frac{1}{2^4}+\frac{1}{3^4}+\frac{1}{4^4}+\frac{1}{5^4}+\textrm{etc.}=\frac{p^4}{90}=Q\\
&1+\frac{1}{2^6}+\frac{1}{3^6}+\frac{1}{4^6}+\frac{1}{5^6}+\textrm{etc.}=\frac{p^6}{945}=R\\
&1+\frac{1}{2^8}+\frac{1}{3^8}+\frac{1}{4^8}+\frac{1}{5^8}+\textrm{etc.}=\frac{p^8}{9450}=S\\
&1+\frac{1}{2^{10}}+\frac{1}{3^{10}}+\frac{1}{4^{10}}+\frac{1}{5^{10}}+\textrm{etc.}=\frac{p^{10}}{93555}=T\\
&1+\frac{1}{2^{12}}+\frac{1}{3^{12}}+\frac{1}{4^{12}}+\frac{1}{5^{12}}+\textrm{etc.}=\frac{691p^{12}}{6825\cdot 93555}=V;
\end{split}
\]
these series are derived from the given rule, with, however, a fair deal of work for the higher powers.
As well, dividing each of these series by the preceding one, we get the following equations:
$p^2=6P=\frac{15Q}{P}=\frac{21R}{2Q}=\frac{10S}{R}=\frac{99T}{10S}=\frac{6825V}{691T}$, etc., where each of these expressions is equal to the square of the perimeter of a circle whose diameter is $1$.

\S 19. Even though the approximate sums of these present series can indeed be worked out easily, they are not a great aid to us in expressing the approximate perimeter of a circle because of the presence of square roots which need to be taken; from the earlier series we have found, we can elicit some expressions which are equal to the perimeter $p$ itself. They are as follows:  
\[
\begin{split}
&p=4\big(1-\frac{1}{3}+\frac{1}{5}-\frac{1}{7}+\frac{1}{9}-\frac{1}{11}+\textrm{etc.}\big)\\
&p=2\Bigg(\frac{1+\frac{1}{3^2}+\frac{1}{5^2}+\frac{1}{7^2}+\frac{1}{9^2}+\frac{1}{11^2}+\textrm{etc.}}{1-\frac{1}{3}+\frac{1}{5}-\frac{1}{7}+\frac{1}{9}-\frac{1}{11}+\textrm{etc.}}\Bigg)\\
&p=4\Bigg(\frac{1-\frac{1}{3^3}+\frac{1}{5^3}-\frac{1}{7^3}+\frac{1}{9^3}-\frac{1}{11^3}+\textrm{etc.}}{1+\frac{1}{3^2}+\frac{1}{5^2}+\frac{1}{7^2}+\frac{1}{9^2}+\frac{1}{11^2}+\textrm{etc.}}\Bigg)\\
&p=3\Bigg(\frac{1+\frac{1}{3^4}+\frac{1}{5^4}+\frac{1}{7^4}+\frac{1}{9^4}+\frac{1}{11^4}+\textrm{etc.}}{1-\frac{1}{3^3}+\frac{1}{5^3}-\frac{1}{7^3}+\frac{1}{9^3}-\frac{1}{11^3}+\textrm{etc.}}\Bigg)\\
&p=\frac{16}{5}\Bigg(\frac{1-\frac{1}{3^5}+\frac{1}{5^5}-\frac{1}{7^5}+\frac{1}{9^5}-\frac{1}{11^5}+\textrm{etc.}}{1+\frac{1}{3^4}+\frac{1}{5^4}+\frac{1}{7^4}+\frac{1}{9^4}+\frac{1}{11^4}+\textrm{etc.}}\Bigg)\\
&p=\frac{25}{8}\Bigg(\frac{1+\frac{1}{3^6}+\frac{1}{5^6}+\frac{1}{7^6}+\frac{1}{9^6}+\frac{1}{11^6}+\textrm{etc.}}{1-\frac{1}{3^5}+\frac{1}{5^5}-\frac{1}{7^5}+\frac{1}{9^5}-\frac{1}{11^5}+\textrm{etc.}}\Bigg)\\
&p=\frac{192}{61}\Bigg(\frac{1-\frac{1}{3^7}+\frac{1}{5^7}-\frac{1}{7^7}+\frac{1}{9^7}-\frac{1}{11^7}+\textrm{etc.}}{1+\frac{1}{3^6}+\frac{1}{5^6}+\frac{1}{7^6}+\frac{1}{9^6}+\frac{1}{11^6}+\textrm{etc.}}\Bigg)
\end{split}
\]

\end{document}